\documentclass[a4paper,11pt,reqno]{amsart}

\usepackage{amsmath}
\usepackage{amsthm}
\usepackage{amsfonts}
\usepackage{amssymb}
\usepackage{mathrsfs}
\usepackage[shortlabels]{enumitem}
\usepackage{graphicx}

\usepackage{mathtools} 
\usepackage{courier} 
\usepackage[T1]{fontenc}  

\usepackage{hyperref}
\usepackage{bm}
\newcommand{\bs}[1] {\bm{#1}}

\newcommand\la\lambda
\newcommand\p\partial
\newcommand\eps\varepsilon
\newcommand{\C}{\mathbb{C}}
\newcommand{\Z}{\mathbb{Z}}
\newcommand\cL{\mathcal L}
\newcommand\cE{\mathcal E}

\newcommand{\sO}{\mathscr{O}}
\renewcommand\d{\,{\rm d}}
\newcommand\e{\mathrm{e}}

\usepackage{amsmath}

\usepackage{cancel}

\makeatletter
\newcommand{\distas}[1]{\mathbin{\overset{#1}{\kern\z@\sim}}}%
\newsavebox{\mybox}\newsavebox{\mysim}
\newcommand{\distras}[1]{%
  \savebox{\mybox}{\hbox{\kern3pt$\scriptstyle#1$\kern3pt}}%
  \savebox{\mysim}{\hbox{$\sim$}}%
  \mathbin{\overset{#1}{\kern\z@\resizebox{\wd\mybox}{\ht\mysim}{$\sim$}}}%
}
\makeatother

\newcommand{\N}{\mathbb{N}}
\newcommand{\R}{\mathbb{R}}

\newcommand{\Hmm}[1]{\leavevmode{\marginpar{\tiny%
$\hbox to 0mm{\hspace*{-0.5mm}$\leftarrow$\hss}%
\vcenter{\vrule depth 0.1mm height 0.1mm width \the\marginparwidth}%
\hbox to
0mm{\hss$\rightarrow$\hspace*{-0.5mm}}$\\\relax\raggedright #1}}}

\newcommand{\cu}{\operatorname{curl}}
\newcommand{\di}{\operatorname{div}}
\newcommand{\gr}{\operatorname{grad}}

\newtheorem{theorem}{Theorem}

\theoremstyle{definition}

\newtheorem{remark}[theorem]{Remark}

\title[]{On a Steklov spectrum in Electromagnetics}
\author[F. Ferraresso]{Francesco Ferraresso$^1$}
\author[P.D. Lamberti]{Pier Domenico Lamberti$^2$}
\author[I. Stratis]{Ioannis G. Stratis$^3$}
\date{\today \\
$^1$School of Mathematics, Cardiff University,  Abacws, Senghennydd Road,
CF24 4AG  Cardiff, United Kingdom, \texttt{FerraressoF@cardiff.ac.uk} \\
$^2$Dipartimento di Tecnica e Gestione dei Sistemi Industriali (DTG),
Universit\`a degli Studi di Padova, Stradella S. Nicola 3,
I-36100 Vicenza, Italy, \texttt{lamberti@math.unipd.it} \\
$^3$Department of Mathematics, National and Kapodistrian University of Athens, Panepistimiopolis,
GR-15784 Zographou, Greece,  \texttt{istratis@math.uoa.gr}
}
\begin{document}

\begin{abstract}
After presenting various concepts and results concerning the classical Steklov eigenproblem, we focus on analogous problems    for time-harmonic Maxwell's equations in a cavity. In this direction we discuss recent rigorous results concerning natural Steklov boundary value problems for the curlcurl operator. Moreover, we explicitly compute eigenvalues and eigenfunctions in the unit ball of the three dimensional Euclidean space by using classical vector spherical harmonics.\\[0.4cm]
\emph{Classification:} 35P15, 35Q61, 78M22. \\[0.1cm]
\end{abstract}

\maketitle

\section{Introduction}

At the macroscopic level, Maxwell's equations read
\begin{equation}
\label{fullmaxwell}
\begin{array}{lll}
\cu \boldsymbol{\widetilde{H}}(x,t) & =  \dfrac{\partial}{\partial t}   \boldsymbol{\widetilde{D}}(x,t) + \boldsymbol{\widetilde{J}}(x,t)&\quad  \mbox{(Amp\`{e}re's law)}, \\
\cu \boldsymbol{\widetilde{E}}(x,t)& =   -\, \dfrac{\partial}{\partial t}  \boldsymbol{\widetilde{B}}(x,t)&\quad \mbox{(Faraday's law of induction)},\\
\di \boldsymbol{\widetilde{D}}(x,t) & =  \widetilde{\varrho}(x,t)&\quad \mbox{(Gauss's law)},\\
\di \boldsymbol{\widetilde{B}}(x,t) & =  0&\quad \mbox{(Gauss's law for magnetism)},
\end{array}
\end{equation}
where $\boldsymbol{\widetilde{E}}$, $\boldsymbol{\widetilde{H}}$ are the electric and the magnetic fields, $\boldsymbol{\widetilde{D}}$, $\boldsymbol{\widetilde{B}}$ are the electric and the magnetic flux densities, $\boldsymbol{\widetilde{J}}$ is the electric current density, and $\varrho$ is the density of the (externally impressed) electric charge.\\
Constitutive relations, i.e., relations of the form $\boldsymbol{\widetilde{D}} =\boldsymbol{\widetilde{D}}(\boldsymbol{\widetilde{E}},\boldsymbol{\widetilde{H}})$, $\boldsymbol{\widetilde{B}} =\boldsymbol{\widetilde{B}}(\boldsymbol{\widetilde{E}},\boldsymbol{\widetilde{H}})$, must accompany \eqref{fullmaxwell}. The constitutive relations for an electromagnetic medium reflect the physics that govern the phenomena
and are expected to comply with the fundamental physical laws, which play the role of physical hypotheses, or postulates,
concerning the properties of the material inside the considered domain. Some constitutive equations are simply phenomenological; others are derived from first principles.
\subsection{Time-harmonic Maxwell's equations in a linear homogeneous isotropic conductive medium}\label{timeharcond}
We now consider electromagnetic wave propagation in a {\it linear homogeneous isotropic} medium in $\mathbb{R}^3$, with  {\em electric permittivity} $\varepsilon$, {\em  magnetic permeability} $\mu$, and {\em electric conductivity} $\sigma$. By our assumptions on the medium, all these three parameters are constant. \\
For such a medium the constitutive relations are
\begin{equation}
 \boldsymbol{\widetilde{D}} = \varepsilon  \boldsymbol{\widetilde{E}}\,,\,\,\,\, \boldsymbol{\widetilde{B}} = \mu  \boldsymbol{\widetilde{H}}.
\end{equation}
The electromagnetic wave with {\it angular frequency} $\,\omega > 0$ will be described by the electric and magnetic field
\begin{equation}
    \boldsymbol{\widetilde{E}}(x,t) = \left(\varepsilon + \dfrac{{\rm i}\,\sigma}{\omega} \right)^{-\,1/2}  {\bf E}(x) \, {\rm e}^{- {\rm i} \omega t},
\end{equation}
\begin{equation}
  \boldsymbol{\widetilde{H}}(x,t) = \mu^{-\,1/2} \,{\bf H}(x) \, {\rm e}^{- {\rm i} \omega t}.
\end{equation}
We additionally consider that the space dependent part ${\bf J}(x)$ of $\boldsymbol{\widetilde{J}}(x,t)$ satisfies
\begin{equation}\label{Ohm}
{\bf J}(x) = \sigma {\bf E}(x)\,\,\,\,\,\,\,\,\,\,\,\,\,\,\,\,\,\,\; \mbox{(Ohm's law)}.
\end{equation}
From Maxwell's equations \eqref{fullmaxwell}, we obtain that the space dependent parts ${\bf E}$ and ${\bf H}$ satisfy the {\it time-harmonic Maxwell equations}
 \begin{equation}\label{thM}
\cu {\bf E}(x) - {\rm i}\,\omega \mu \, {\bf H}(x) = {\bf 0}\,, \,\,\,\,
\cu {\bf H}(x) + {\rm i}\,\omega\varepsilon\, {\bf E}(x) = {\bf 0}.
\end{equation}
Let us note that since $\varepsilon$ and $\mu$ are constant we have  that both the fields ${\bf E}$ and ${\bf H}$ are divergence-free:
\begin{equation}\label{divfree}
 \di {\bf E} = \di {\bf H} = 0.
\end{equation}

Let $\Omega$ be a bounded domain (i.e. a bounded connected open set) in $\R^3$ with smooth boundary $\Gamma$.
We consider the following classical  boundary value problem that
involves the  ``perfect conductor'' condition on $\Gamma$
\begin{equation}\label{PEC}
\left\{
\begin{array}{ll}
\cu {\bf E} - {\rm i} \omega \mu {\bf H} = {\bf 0}\,,\,\, \cu {\bf H} + {\rm i} \omega \varepsilon {\bf E} = {\bf 0} ,& \ \ {\rm in}\ \Omega ,\\
{\bm \nu} \times {\bf E} = {\bf m} ,& \ \ {\rm on}\ \Gamma ,
\end{array}
\right.
\end{equation}
where ${\bm \nu}$ denotes the unit outer normal to $\Gamma$.  By eliminating ${\bf H}$ we obtain
\begin{equation}\label{PECH}
\left\{
\begin{array}{ll}
\cu\cu {\bf E} - k^2{\bf E} = {\bf 0} , & \ \ {\rm in}\ \Omega ,\\
{\bm \nu} \times {\bf E} = {\bf  m} ,& \ \ {\rm on}\ \Gamma ,
\end{array}
\right.
\end{equation}
where $k^2:= \omega^2\varepsilon\mu $, and as usual we assume that  ${\rm Im} k\geq 0$.

Instead of the standard interior Calder\'{o}n operator (see e.g., \cite{ces}, \cite{kswy}), in what follows we consider its variant  defined by ${\bf m}\mapsto ({\bm \nu }\times {\bf H} ) \times {\bm \nu }$, i.e.
\begin{equation}
{\bm \nu} \times {\bf E}\mapsto  -\frac{{\rm i}}{\omega \mu}  ({\bm \nu} \times \cu {\bf E}) \times {\bm \nu} \, .
\end{equation}
  Calder\'{o}n operators are also called Poincar\'{e}-Steklov, or impedance, or admittance, or capacity operators.

\subsection{The classical Steklov eigenvalue problem}

We recall that the classical Steklov\footnote{Vladimir Andreevich Steklov (1864 - 1926) was not only an outstanding
mathematician, who made many important contributions to Applied
Mathematics, but also had an unusually bright personality.
The Mathematical Institute of the Russian Academy of Sciences in Moscow
bears his name.
On his life and work see the very interesting paper \cite{kukukwnapoposi}.
} eigenvalue problem on a bounded domain $\Omega$ of $\R^n$, $n\geq 2$,  is the problem
\begin{equation}
\label{classicsteklov}
\left\{\begin{array}{ll}
\Delta u =0,&\  {\rm in }\ \Omega,\\
\dfrac{\partial u}{\partial \nu}  =\lambda u,&\  {\rm on }\ \partial\Omega
\end{array}
\right.
\end{equation}
in the unknowns $u$ (the eigenfunction) and $\lambda$ (the eigenvalue). The domain $\Omega$ is assumed to be sufficiently regular (usually one requires that the boundary $\Gamma$ is at least Lipschitz continuous) and the (scalar)  harmonic function  $u$ is required to belong to the standard Sobolev space $H^1(\Omega)$. This problem can be considered as the eigenvalue problem for the celebrated Dirichlet-to-Neumann map defined as follows. Given the solution  $u\in H^1(\Omega)$ to the  the Dirichlet problem\
\begin{equation}
\left\{\begin{array}{ll}
\Delta u =0,&\  {\rm in }\ \Omega,\\
u=f,&\  {\rm on }\ \partial\Omega
\end{array}
\right.
\end{equation}
with datum $f\in H^{1/2}(\Gamma)$,  one can consider the normal derivative $\frac{\partial u}{\partial \nu}$ of $u$ as an element of $H^{-1/2}(\Gamma)$,
where   $H^{1/2}(\Gamma)$ is the standard Sobolev space defined on $\Gamma$  and $H^{-1/2}(\Gamma)$ its dual. This allows to define the map ${\mathcal{D}}$ from $H^{1/2}(\Gamma)$  to $H^{-1/2}(\Gamma)$ by setting
$$
{\mathcal{D}} f = \dfrac{\partial u}{\partial \nu}\,.
$$
The map ${\mathcal{D}}$ is called Dirichlet-to-Neumann map and its eigenpairs $(f, \lambda )$ correspond to the eigenpairs $(u, \lambda)$ of problem \eqref{classicsteklov}, $f$ being the trace of $u$ on $\Gamma$.

\subsection{The electromagnetic Steklov eigenvalue problem}

The natural analogue in electromagnetics of the classical Steklov   problem \eqref{classicsteklov}   in $\R^3$  can be defined as the eigenvalue problem for the (rescaled) interior Calder\'{o}n operator defined above, namely the map
\begin{equation}
{\bm \nu} \times {\bf E}\mapsto    - ({\bm \nu} \times \cu {\bf E}  ) \times {\bm\nu }  .
\end{equation}
Therefore, one looks for values $\lambda$ such that
$
 ({\bm \nu} \times \cu {\bf E} ) \times {\bm \nu} = - \lambda {\bm \nu} \times {\bf E},
$
or, equivalently (by taking another cross product by $\bm \nu $)
\begin{equation}
 {\bm \nu} \times \cu {\bf E}  =  \lambda  {\bf E}_{\rm T},
 \end{equation}
where $ {\bf E}$ satisfies the equation      $\cu\cu {\bf E} -  k^2 {\bf E} = {\bf 0}  $   and $ {\bf E}_{\rm T}    :=   ({\bm \nu}\times  {\bf E} ) \times   {\bm \nu}    $ is the tangential component of ${\bf E}$.  In conclusion, the Steklov eigenvalue problem for Maxwell's equations is
\begin{equation}\label{stemax}
\left\{
\begin{array}{ll}
\cu\cu {\bf E} -  k^2 {\bf E} = {\bf 0} , & \ \ {\rm in}\ \Omega ,\\
 {\bm \nu} \times \cu {\bf E}  =  \lambda  {\bf E}_{\rm T} ,& \ \ {\rm on}\ \Gamma .
\end{array}
\right.
\end{equation}

To the best of our knowledge problem \eqref{stemax} was first introduced for $k>0$  by  J. Caman\~o, C. Lackner and P.  Monk in \cite{calamo} where it was pointed out that the spectrum of this problem is not discrete. In particular, for the case of the unit ball in $\R^3$ it turns out that the eigenvalues consist of two infinite sequences, one of which is divergent and the other is converging to zero. To overcome this issue, in that paper a modified problem, having discrete spectrum,  is considered and then used to study an inverse scattering problem.

On the other hand, in \cite{lamstr}, two of the authors of the present paper have analyzed problem  \eqref{stemax} only for tangential vector fields ${\bf E}$ in which case the problem can be written in the form
\begin{equation}\label{stemaxtan}
\left\{
\begin{array}{ll}
\cu\cu {\bf E} -  k^2 {\bf E} = {\bf 0} , & \ \ {\rm in}\ \Omega ,\\
 {\bm \nu} \times \cu {\bf E}  =  \lambda  {\bf E} ,& \ \ {\rm on}\ \Gamma .
\end{array}
\right.
\end{equation}

Note that the boundary condition in \eqref{stemaxtan} automatically implies that  ${\bf E}$ is tangential, that is $ ${\bf E}$\cdot {\bm \nu}=0$ on $\Gamma$. Because of this restriction, the null sequence of eigenvalues disappears and the spectrum turns out to be discrete.

Steklov eigenproblems  have been and are extensively studied for a variety of differential operators, linear and nonlinear, mostly in the scalar case. On the contrary, there  are not so many publications devoted to analogous problems for Maxwell's equations:  apart the two papers mentioned above, we are aware of the ones by F. Cakoni, S. Cogar and P. Monk \cite{cacomon},  by S. Cogar \cite{cogar2}, \cite{cogar3}, \cite{cogar4}, by S. Cogar, D. Colton and P. Monk \cite{cocomo2019}, by S. Cogar,  P. Monk \cite{como2020}, and by M. Halla \cite{halla1}, \cite{halla2}.\\

In this paper, after presenting several comments and results for the classical Steklov eigenproblem (see Section~\ref{classicalsec}),  we discuss the approach introduced in  \cite{lamstr} (see Section~\ref{lamstrsec})  and we include explicit computations of the eigenvalues of the related eigenvalue problems in the case of the unit ball in $\R^3$ (see Section~\ref{ballsec}).

\section{On the classical Steklov eigenvalue problem}
\label{classicalsec}

\subsection{Some indicative applications}

 Problem \eqref{classicsteklov} has a long history which goes back to the paper \cite{st} written by Steklov himself. It is not easy to give a complete account of all possible applications of problem and its variants. Here we briefly mention three main fields of investigation.

 \begin{itemize}
 \item {\it The sloshing problem.}  It consists in the study of small oscillations of a liquid in finite basin which can be thought as bounded container (a tank, a mug, a snifter etc.). The basin is represented by a bounded domain $\Omega$ in ${\mathbb{R}}^3$ with boundary $\Gamma =\Gamma_1\cup \Gamma_2$, where $\Gamma_1$ is a two dimensional domain representing the horizontal (free) surface of the liquid at rest, and $\Gamma_2$ represents the bottom of the basin. In this case the Steklov boundary condition $\frac{\partial u  }{\partial \nu}=\lambda u $ is imposed only on $\Gamma_1$, while the Neumann condition $\frac{\partial u  }{\partial \nu}$ is imposed on $\Gamma_2$.   The gradients $\gr u (x,y,z)$ of solutions $u$ represent the (stationary) velocity fields of the oscillations and $\sqrt {\lambda }$ the corresponding frequency. In particular $u(x,y,0) $ is proportional to the elevation of the free surface and the so-called ``high-spots''' correspond to its maxima.
 We refer to  \cite{canavati},   \cite{chiado}, \cite{kukukwnapoposi}, \cite{levitin} for more details on classical and more recent aspects of the problem.

\item {\it Electrical prospection}  The study of the Dirichlet-to-Neumann map received a big impulse from the seminal paper \cite{cal} by Calder\'{o}n which poses the inverse problem of recovering the electric conductivity $\gamma$ of an electric body $\Omega$ from the knowledge of (the energy form associated with)
the voltage-to-current map $D_{\gamma}$ defined in the same way as
$$
D_{\gamma} f = \gamma     \dfrac{\partial u_{\gamma}}{\partial \nu},
$$
where $u_{\gamma}$ is the solution to the problem
\begin{equation}
\label{classicsteklovgamma}
\left\{\begin{array}{ll}
{\rm div} ( \gamma \gr  u_{\gamma} )=0,&\  {\rm in }\ \Omega,\vspace{2mm}\\
  u_{\gamma}    =f,&\  {\rm on }\ \Gamma  .
\end{array}
\right.
\end{equation}
The problem of Calder\'{o}n was solved in the fundamental paper \cite{syuh}, where it is proved that the map 
\[f\mapsto \int_{\Gamma} f  D_{\gamma} f\,d\sigma \,= \int_{\Omega}\gamma |\gr u_{\gamma}|^2\,dx\]
uniquely identifies $\gamma$ (in the class of conductivities $\gamma$ of class $C^{\infty}(\overline\Omega)$).

\item{\it Vibrating membranes}   Problem \eqref{classicsteklov} can be used in linear elasticity to model the vibrations of a free membrane $\Omega$ in $ {\mathbb{R}}^2$ with mass concentrated at the boundary. Recall that the normal modes of a free membrane with mass density $\rho$ are the solutions to the Neumann eigenvalue problem
\begin{equation}\label{neustrip}
\left\{ \begin{array}{ll}
-\Delta u= \lambda \rho u,&\  {\rm in }\ \Omega,\\
  \dfrac{\partial u}{\partial \nu}   =0,&\  {\rm on }\ \Gamma .
\end{array}
\right.
\end{equation}
The total mass of the membrane is given by $M=\int_{\Omega} \rho (x)\,dx$. If we consider a family of mass densities $\rho_{\epsilon}$ for $\epsilon>0$, such
that the support of $\rho_{\epsilon}$ is contained in a neighborhood of the boundary  $\Gamma$ of radius $\epsilon$ (with $\rho_{\epsilon}$ constant therein) and such that the total mass $M_{\epsilon}=M$ does not depend on $\epsilon$, then the solutions of problem \eqref{neustrip} converge to the solutions of
\begin{equation}
\label{classicsteklovrho}
\left\{\begin{array}{ll}
\Delta u =0,&\  {\rm in }\ \Omega,\\
 \dfrac{\partial u}{\partial \nu}  =   \lambda \rho u,&\  {\rm on }\ \Gamma ,
\end{array}
\right.
\end{equation}
where $\rho = M/  |\Gamma |$ and $|\Gamma|$ is the perimeter of $\Gamma$. Thus problem \eqref{classicsteklovrho} can be considered as a limiting/critical case of a family of Neumann eigenvalue problems. We refer to \cite{dalla, lampro1, lampro2} for further details, in particular for an asymptotic analysis.
 \end{itemize}

 Finally, we mention that the Steklov problem has been recently used in \cite{hey} for a mathematical model related to  the study of information transmission in the neural network of the human brain.

\subsubsection{Details on the formulation and its connections to trace theory}
\label{tracesec}

For any $n\geq 2$ the weak (variational) formulation of problem \eqref{classicsteklov} is easily obtained by multiplying the equation $\Delta u=0 $ by a test function $\varphi $ and integrating by parts over $\Omega$. This simple computation leads to the equality
\begin{equation}
\label{weakclassic}
\int_{\Omega}\gr u\cdot \gr \varphi\, dx= \lambda \int_{\Gamma} u\varphi\, d\sigma
\end{equation}
which needs to be satisfied for all functions  $\varphi \in H^1(\Omega)$, and can be taken as the formulation of the classical problem \eqref{classicsteklov} in $H^1(\Omega)$. The advantage of  formulation \eqref{weakclassic} is evident since it easily allows to apply standard  tools from functional analysis and calculus of variations to prove existence of solutions.
Indeed, one can prove that the spectrum of problem \eqref{weakclassic} is discrete and consists of a divergent sequence of eigenvalues
$$
0=\lambda_1\le \lambda_2\le \dots \le \lambda_k \le \dots
$$
where it is assumed that each eigenvalue is repeated as many times as its multiplicity, which is finite. The corresponding eigenfunctions $\varphi_j$ define a complete orthogonal system for the subspace $\mathsf{H} (\Omega )$ of harmonic functions in $H^1(\Omega)$. Note that  $\mathsf{H} (\Omega )$ can be described as
$\mathsf{H} (\Omega )=\{u\in H^1(\Omega):\ \int_{\Omega}\gr u\cdot \gr \varphi \,dx =0,\ \forall \varphi \in C^{\infty}_c(\Omega) \}$. In fact, the following decomposition holds
\begin{equation}
\label{decomposition}
H^1(\Omega) =H^1_0(\Omega)\oplus \mathsf{H} (\Omega )\, .
\end{equation}

We describe one straightforward way to prove these results since this method will be applied to the case of Maxwell's equations.

By adding the term $ \int_{\Gamma}u\varphi \,d\sigma$ to both sides of the equation \eqref{weakclassic}
and setting $\mu =\lambda +1$, one gets the equation
$$
\int_{\Omega}\gr u\cdot \gr \varphi \,dx + \int_{\Gamma}u\varphi \,d\sigma= \mu \int_{\Gamma} u\varphi \,d\sigma
$$
where the quadratic form associated to the left-hand side, namely
$$
Q(u):=\int_{\Omega}|\gr u|^2 dx + \int_{\Gamma}|u|^2 \,d\sigma
$$
is coercive in $H^1(\Omega )$ and in particular defines a norm $Q(u)^{1/2}$ equivalent to the Sobolev norm of $H^1(\Omega)$. Thus, the operator
${\mathcal{L}}$ from $H^1(\Omega)$ to its dual defined by the pairing
$$
\langle{\mathcal{L}}u, \varphi \rangle = \int_{\Omega}\gr u\cdot \gr \varphi \,dx + \int_{\Gamma}u\varphi \,d\sigma
$$
 is invertible. Then we can consider the operator $T$ from $H^1(\Omega)$ to itself defined by $T={\mathcal{L}}^{-1}\circ J$ where
 $J$ is the operator from $H^1(\Omega)$ to its dual defined by the pairing
 $$
 \langle Ju, \varphi \rangle = \int_{\Gamma}u\varphi \,d\sigma .
 $$
It is simple to see that $T$ is a selfadjoint operator with respect to the scalar product associated with the quadratic form $Q$ above (see e.g., \cite{la2014} for an analogous Steklov problem). Moreover,  since the standard trace operator $Tr$  from   $H^1(\Omega)$ to $L^2(\Gamma )$ is compact, it follows that $T$ is also compact, hence its spectrum consists of zero and a decreasing divergent sequence of positive eigenvalues $\mu_j$.
Thus the eigenvalues $\lambda_j$ above can be defined by the equality $\mu_j = (\lambda_j+1)^{-1}$. Moreover, since the kernel of $T$ is exactly $H^1_0(\Omega)$ and its orthogonal is $\mathsf{H} (\Omega )$, the decomposition \eqref{decomposition} immediately follows.

If $\Omega$ is the ball of radius $R$ centred at zero, the eigenvalues are given by all numbers of the form
$$
l_j=\frac{j}{R},\ \ j\in \N_0.
$$
The corresponding eigenfunctions are the homogeneous polynomials of degree $j$ and can be written in spherical coordinates in  the form
$$
u(r, \xi ) = r^jY_j(\xi )
$$
for $r=|x|\geq 0$ and $\xi =x/ |x|\in S^{n-1}$ (the $(n-1)$-dimensional unit sphere), where  $Y_j$ is any spherical harmonic of degree $j$. In particular,
 the multiplicity of $l_j$ is
$(2j+n-2) (j+n-3)! /(j!(n-2)!)$,
and only $l_0$ is simple, the corresponding eigenfunctions being the constant functions.   Note that the enumeration  $l_j$, $j\in \N_0$ is different from the enumeration $\lambda_j$, $j\in \N$ discussed above since it does not take into account the  multiplicity of the eigenvalues. We note {\it en passant} that the eigenvalues of the Laplace-Beltrami operator on the $(n-1)$-dimensional sphere of radius $R$ are given by the formula $\sigma_l= l(l+n-2)/R^2$ and coincide with the squares of the $l_j$ for $n=2$ (the corresponding  eigenfunctions are given by  the restrictions of the corresponding Steklov eigenfunctions). We refer to \cite{pro} for further discussions.

We would like now to describe the method of Auchmuty~\cite{Auch2006} for the spectral representation of the trace  space $H^{1/2}(\Gamma)$, since the same method will be used in the vectorial case (in which case the Steklov eigenvectors for  Maxwell's equations will be used).

By exploiting  an argument similar to the one discussed above (with an operator analogous to $T$  defined on $L^2(\Gamma)$ rather than on $H^1(\Omega)$), one can actually see that the traces of the  eigenfunctions $\varphi_j$ on $\Gamma$ define a complete orthogonal system for $L^2(\Gamma)$.
Here for simplicity we also write $\varphi_j$ instead of $Tr (\varphi_j)$. Assume that those eigenfunctions are normalized in $L^2(\Gamma)$, that is $\int_{\Gamma}|\varphi_j|^2\,d\sigma=1$. Then,  by equation \eqref{weakclassic} it follows that $\varphi_j/\sqrt{\lambda_j+1}$ is normalized in $H^1(\Omega)$, that is $Q(\varphi_j/\sqrt{\lambda_j+1})=1$. Thus, $\mathsf{H}(\Omega)$ can be described as follows
\begin{equation}
\label{expansion}
\mathsf{H}(\Omega)=\left\{  \sum_{j=1}^{\infty} c_j \frac{\varphi_j}{\sqrt{\lambda_j+1}}: \ \sum_{j=1}^{\infty}|c_j|^2< \infty  \right\} .
\end{equation}
Recall that the  trace space  $Tr (H^1(\Omega))$   coincides with the standard fractional Sobolev space $H^{1/2}(\Gamma)$ and note that $Tr (H^1(\Omega))=  Tr ( \mathsf{H}(\Omega))$  by \eqref{decomposition}.  This, combined with \eqref{expansion}, yields
\begin{equation}\label{auchmuty}
\begin{split}
H^{1/2}(\Gamma) =
Tr ( \mathsf{H}(\Omega)) &=\left\{  \sum_{j=1}^{\infty} c_j \frac{Tr (\varphi _j)}{\sqrt{\lambda_j+1}}: \ \sum_{j=1}^{\infty}|c_j|^2< \infty  \right\} \\
&= \left\{  \sum_{j=1}^{\infty} c_j Tr( \varphi _j): \ \sum_{j=1}^{\infty}(\lambda_j+1)|c_j|^2< \infty  \right\}.
\end{split}
\end{equation}

If $\Omega$ is sufficiently smooth then Weyl's law, describing  the asymptotic behavior of the eigenvalues, is
$$
\lambda_j\sim c j^{\frac{1}{n-1}},\ \ {\rm as }\  j\to \infty
$$
where $c$ is an explicitly known constant. It follows that  the space $H^{1/2}(\Gamma) $ can be described as folllows
$$
H^{1/2}(\Gamma) =
\left\{  \sum_{j=1}^{\infty} c_j  \varphi _j: \ \sum_{j=1}^{\infty}j^{\frac{1}{n-1}}|c_j|^2< \infty  \right\}.
$$
Thus the condition on the Fourier coefficients $c_j$ is that the sequence
\begin{equation}
\label{weighted}
j^{\frac{1}{2(n-1)}}c_j,\ \ j\in {\mathbb{N}}
\end{equation}
belongs to the space $\ell^2$ of square summable sequences.  We note that  the  appearance  of the factor  $1/2$ at  the exponent in \eqref{weighted}  is not artificial and corresponds to the exponent of the space  $H^{1/2}(\Gamma)$. In fact  an analogous representation was found in \cite{lampro3} for the space $H^{3/2}(\Gamma)$ where the exponent $3/2$ naturally appears by using the Weyl asymptotic for a  biharmonic Steklov eigenvalue  problem.

\section{On the electromagnetic Steklov eigenproblem}
\label{lamstrsec}

In this section we briefly present some of the results in \cite{lamstr}  for problem \eqref{stemaxtan}. In the sequel $\Omega$ will denote
a bounded domain in $\R^3$ with sufficiently smooth boundary, say of class $C^{1,1}$ (see e.g., \cite[Definition 1]{lamzac}).  As done in \cite{coda} for analogous problems, we introduce a  penalty term $\theta\, \gr\,  \di u$  in the equation, where $\theta$ can be any positive number,  in order to guarantee the coercivity of the quadratic form associated with the corresponding differential operator.  Namely, we consider the eigenvalue problem

\begin{equation}\label{classiccucu1}
\left\{
\begin{array}{ll}
\cu\, \cu \bm{E} -k^2 \bm{E} -\theta\, \gr\,  \di \bm{E}=\bm{0},& \ \ {\rm in}\ \Omega ,\\
\bm{\nu} \times \cu \bm{E} = \lambda \bm{E},& \ \ {\rm on}\ \Gamma ,\\
\bm{E} \cdot \bm{\nu}=0, & \ \ {\rm on}\ \Gamma
\end{array}
\right.
\end{equation}
where $\bm{E}$ is the unknown vector field. Here we allow $k^2\in \R$ to be not necessarily positive.  Recall that the second boundary condition above is in fact embodied in the first one but we prefer to emphasize it since we need to include it in the definition of the energy space.\\

 By $L^2(\Omega)$, $H^1(\Omega)$, $H^{1}_{0}(\Omega)$, $L^2(\Gamma)$, $H^{1/2}(\Gamma), H^{-1/2}(\Gamma)$,
we denote the standard Lebesgue and Sobolev spaces. We also employ the following spaces:

\begin{itemize}
\item $H(\cu, \Omega) = \{{\bm u} \in (L^2(\Omega))^3 : \cu {\bm u} \in (L^2(\Omega))^3 \}\,,$\\
with norm: $ \|{\bm u}\|_{H(\cu, \Omega)} = \left( \|{\bm u}\|^{2}_{(L^2(\Omega))^3} + \|\cu {\bm u}\|^{2}_{(L^2(\Omega))^3} \right)^{1/2} $
\item $H(\di, \Omega) = \{{\bm u} \in (L^2(\Omega))^3 : \di {\bm u} \in L^2(\Omega) \}\,,$\\
with norm: $ \|{\bm u}\|_{H(\di, \Omega)} = \left( \|{\bm u}\|^{2}_{(L^2(\Omega))^3} + \|\di {\bm u}\|^{2}_{L^2(\Omega)} \right)^{1/2} $
\item $H_{0}(\di, \Omega) = \{{\bm u} \in H(\di, \Omega): {\bm \nu} \cdot {\bm u} = 0 \,\, {\rm on}\ \Gamma \}$
\item $X_{\rm \scriptscriptstyle T}(\Omega) = H(\cu, \Omega) \cap H_0(\di, \Omega)\,,$  with norm:\\
$ \|{\bm u}\|_{X_{\rm \scriptscriptstyle T}(\Omega)} = \left( \|{\bm u}\|^{2}_{(L^2(\Omega))^3} + \|\cu {\bm u}\|^{2}_{(L^2(\Omega))^3} + \|\di {\bm u}\|^{2}_{L^2(\Omega)} \right)^{1/2} $
\end{itemize}

 We refer to \cite{acl}, \cite{ces}, \cite{dali3}, \cite{gira}, \cite{kihe}, \cite{monk}, \cite{rsy} for details.

 It is important to note that since we have assumed  $\Omega$ to be of class $C^{1,1}$, the space $X_{\rm \scriptscriptstyle T}( \Omega)$ is continuously embedded in $(H^1(\Omega ))^3$ and there exists
$c>0$ such that the {\it Gaffney inequality}
$$\| {\bm u} \|_{(H^1(\Omega))^3}\le c \, \left(    \|   {\bm u}\|_{L^2(\Omega)^3 } +     \|  \cu {\bm u}\|_{L^2(\Omega)^3} +  \|  \di {\bm u}\|_{L^2(\Omega)} \right),$$
holds  for all ${\bm u}\in   X_{\rm \scriptscriptstyle T}( \Omega)$.\\

 Problem \eqref{classiccucu1} has to be interpreted in the weak sense as follows:
find $\bm{E}\in X_{\rm \scriptscriptstyle T}(\Omega)$  such that
\begin{equation}\label{weakcucu1e}
\int_{\Omega}\cu \bm{E}\cdot  \cu \bm{\varphi}\, dx -k^2 \int_{\Omega}\bm{E} \cdot  \bm{\varphi}\, dx + \theta \int_{\Omega }\di \bm{E}\,  \di  \bm{\varphi} \, dx =  -  \lambda \, \int_{\Gamma} \bm{E} \cdot  \bm{\varphi} \, d\sigma\, ,
\end{equation}
for all $ \bm{\varphi} \in X_{\rm \scriptscriptstyle T}(\Omega)$.

The above  formulation is obtained from \eqref{classiccucu1}  by a standard procedure: for a smooth solution $\bm{E}$ of  \eqref{classiccucu1},  we multiply both sides of the first equation in \eqref{classiccucu1} by $\bm{\varphi} \in X_{\rm \scriptscriptstyle T}(\Omega)$, integrate by parts and use the following standard Green-type formula
\begin{equation}\label{parts}
\int_{\Omega}\cu \bm{E}\cdot  \cu \bm{\varphi}\, dx =\int_{\Omega}\cu\, \cu \bm{E} \cdot \bm{\varphi} \,dx-\int_{\Gamma }(\nu \times \cu \bm{E} )\cdot \bm{\varphi}\, d\sigma\, .
\end{equation}
Conversely, by the Fundamental Lemma of the Calculus of Variations (see e.g., \cite{brezis}), one can see that if   $\bm{E}$ is a smooth solution of $\eqref{weakcucu1e}$ then it is also a solution of $\eqref{classiccucu1}$ in the classical sense.

Note that the weak formulation  allows to avoid assuming additional regularity assumptions on $\Gamma$, see  e.g., \cite{weber}.

In order to study our eigenvalue problem, we need to assume that $k^2$ does not coincide with  an eigenvalue $A$ of the problem

\begin{equation}\label{dirichlet}
\left\{
\begin{array}{ll}
\cu\, \cu \bm{E} -\theta\, \gr\,  \di \bm{E}=  A\bm{E} ,& \ \ {\rm in}\ \Omega ,\\
\bm{\nu} \times  \bm{E} = \bm{0},& \ \ {\rm on}\ \Gamma ,\\
\bm{E} \cdot \bm{\nu}=0, & \ \ {\rm on}\ \Gamma .
\end{array}
\right.
\end{equation}
Clearly the two boundary conditions above are equivalent to the Dirichlet condition $\bm{E}=\bm{0}  $ on $\Gamma$.

We note that \eqref{dirichlet} has a discrete spectrum which consists of a sequence $A_n$, $n\in \N$ of positive eigenvalues of finite multiplicity, the first one being
 \begin{equation}
\label{ray}
A_1=\min_{\substack{\bm{\varphi}\in (H^1_0(\Omega))^3\  \bm{\varphi}\ne \bm{0} }  }\frac{  \int_{\Omega}|\cu \bm{\varphi}|^2\, dx + \theta \, \int_{\Omega }|\di \bm{\varphi}|^2 \, dx   }{\int_{\Omega} |\bm{\varphi}|^2\,dx}  >0.
\end{equation}

For the sake of brevity, we assume in the sequel that $k^2<A_1$.
For details on the more general case
\begin{equation}
\label{exclusion}
A_n<k^2<A_{n+1}
\end{equation}
we refer to \cite{lamstr}.\\

 The key result in the considered case is the following.

  \begin{theorem}\label{mainspectrum}   Let $k^2 <A_1$ and $\theta >0$.
  The  eigenvalues of problem \eqref{classiccucu1} are real, have finite multiplicity and can be represented by a sequence $\lambda_n,\ n \in \N$, divergent to $-\infty$. Moreover, the following min-max representation holds:
  \begin{equation}\label{minmax}
  \lambda_n=  - \min_{ \substack{ V\subset X_{\rm \scriptscriptstyle T}(\Omega )  \\ {\rm dim }V=n }  }\  \,  \max _{\bm{\varphi}\in V\setminus (H^1_0(\Omega ))^3}
  \frac{  \int_{\Omega} \left( |\cu \bm{\varphi}|^2 -k^2  |\bm{\varphi}|^2 + \theta  |\di \bm{\varphi} |^2 \right) dx }{\int_{\Gamma} | \bm{\varphi}  |^2\, dx}\, .
  \end{equation}
  \end{theorem}

To prove this result we follow the strategy described in Section~\ref{tracesec}. Namely,  by adding the term $\eta \int_{\Gamma}  \bm{E}\cdot \bm{\varphi} \, d\sigma$ to  both sides of equation \eqref{weakcucu1e} we obtain

\begin{eqnarray}\label{weakcucu1eta}  \lefteqn{
\int_{\Omega} \! \cu \bm{E} \cdot \cu \bm{\varphi}\, dx \, -\, k^2\!  \int_{\Omega}\! \bm{E} \cdot \bm{\varphi}\, dx \, + \, \theta\!  \int_{\Omega }\! \di \bm{E}\,  \di \bm{\varphi} \, dx }  \nonumber  \\
& & \qquad\qquad  \qquad\qquad  \qquad\qquad
+\, \eta \int_{\Gamma}\! \bm{E}\cdot \bm{\varphi} \, d\sigma =
 \gamma\!  \int_{\Gamma}\! \bm{E} \cdot  \bm{\varphi} \, d\sigma
\end{eqnarray}
where $\gamma = -\lambda +\eta$. Under our assumptions,  it is proved in \cite[Thm.~3.1]{lamstr} that if $\eta $ is big enough then the quadratic form associated with the left-hand side of equation \eqref{weakcucu1eta}, that is
$$
{\mathcal{Q}}(\bm{E}):= \int_{\Omega}\!  |\cu \bm{E} |^2\, dx \, -\, k^2\!  \int_{\Omega}\! |\bm{E} |^2\, dx \, + \, \theta\!  \int_{\Omega } \! |\di \bm{E}|^2 \, dx \, +\, \eta\! \int_{\Gamma} |\bm{E} |^2 \, d\sigma ,
$$
is coercive in $X_{\rm \scriptscriptstyle T}(\Omega )$, hence $({\mathcal{Q}}(\bm{E}))^{1/2}$ defines a norm equivalent to that of $X_{\rm \scriptscriptstyle T}(\Omega )$.

Thus, the operator
${\mathcal{L}}^{\eta}$ from $X_{\rm \scriptscriptstyle T}(\Omega )$ to its dual defined by the pairing
\begin{multline*}
\langle  {\mathcal{L}}^{\eta}\bm{E}, \bm{\varphi} \rangle := \int_{\Omega}\!\cu \bm{E} \cdot \cu \bm{\varphi}\, dx  - k^2\! \! \int_{\Omega}\!\bm{E} \cdot \bm{\varphi}\, dx  \\
+  \theta\!  \int_{\Omega }\!\di \bm{E}\,  \di \bm{\varphi} \, dx  + \eta\! \int_{\Gamma} \! \bm{E}\cdot \bm{\varphi} \, d\sigma
\end{multline*}
is invertible.

Then we can consider the operator ${\mathcal{T}}$ from $X_{\rm \scriptscriptstyle T}(\Omega )$ to itself, defined by
$${\mathcal{T}}=({\mathcal{L}}^{\eta})^{-1}\circ \mathcal{J}$$
where
 $ \mathcal{J}$ is the operator from $X_{\rm \scriptscriptstyle T}(\Omega )$ to its dual defined by the pairing
 $$
 \langle  \mathcal{J} \bm{E}, \bm{\varphi }\rangle = \int_{\Gamma}  \bm{E}\cdot \bm{\varphi } \,d\sigma
 $$
 for all $ \bm{E}, \bm{\varphi }\in    X_{\rm \scriptscriptstyle T}(\Omega )$. As in the case described  in Section~\ref{tracesec}, it is not difficult to prove that  ${\mathcal{T}}$ is a selfadjoint operator with respect to the scalar product associated with the quadratic form ${\mathcal{Q}}$ above. Again,   since the trace operator is compact, it follows that ${\mathcal{T}}$ is also compact, hence its spectrum consists of zero and a decreasing divergent sequence of positive eigenvalues $\gamma_j$.
Thus the eigenvalues $\gamma_j$ above can be defined by the equality $\gamma_j = (-\lambda_j+\eta)^{-1}$. Then the characterization  in \eqref{minmax} follows by the classical Min-Max Principle applied to the operator  ${\mathcal{T}}$.

It follows by the previous results that the space $X_{\rm \scriptscriptstyle T}(\Omega)$ can be decomposed as an orthogonal sum with respect to the scalar product associated with the form $\mathcal{Q}$,  namely
$$
X_{\rm \scriptscriptstyle T}(\Omega )=  {\rm Ker}  {\mathcal{T}}   \oplus   \left( {\rm Ker}  {\mathcal{T}}  \right)^{\perp } = (H^1_0(\Omega ))^3\oplus   \mathcal{H}(\Omega )\, .
$$
where

\begin{equation}\label{armoniche}
\begin{split}
&\mathcal{H}(\Omega ):= \left( {\rm Ker}  {\mathcal{T}}\right)^{\perp}   =\left\{ \bm{E}\in X_{\rm \scriptscriptstyle T}(\Omega ):\
\int_{\Omega}\cu  \bm{E} \cdot \cu \bm{\varphi}\, dx   \right. \\
&\qquad \quad\left.
-k^2\int_{\Omega}\bm{E}\cdot \bm{\varphi} \, dx+ \theta \int_{\Omega }\di \bm{E}\,  \di \bm{\varphi} \, dx= 0,\  \forall  \bm{\varphi} \in   (H^1_0(\Omega ))^3   \right\}.
\end{split}
\end{equation}

Note that  $\bm{E}\in  {\mathcal{H}}(\Omega ) $
 if and only if $\bm{E}$ is a weak solution in $(H^1(\Omega))^3$ of the problem

\begin{equation}\label{weakharmonic}
\left\{
\begin{array}{ll}
\cu\, \cu \bm{E} -k^2 \bm{E} -\theta\, \gr\,  \di \bm{E}=\bm{0},& \ \ {\rm in}\ \Omega ,\\
\nu \cdot \bm{E} = 0 ,& \ \ {\rm on}\ \Gamma .
\end{array}
\right.
\end{equation}

Solutions to problem \eqref{weakharmonic} play  the same role as the harmonic functions in $\mathsf{H}(\Omega)$ used in Section~\ref{tracesec} and, similarly,  the eigenfunctions associated with the eigenvalues $\gamma_n$ define a complete orthonormal system of $ {\mathcal{H}}(\Omega )$.\\

Let $\Sigma =\{ \lambda_n:\ n\in \N\}$.  It is important to known whether $0\in \Sigma$. This condition can be clarified as follows.
We consider two auxiliary eigenproblems.  The first is the classical eigenvalue problem for the Neumann Laplacian
\begin{equation}\label{neulap}
\left\{
\begin{array}{ll}
-\Delta \phi =\lambda \phi,&\ \ {\rm in }\, \Omega,\vspace{1mm}\\
\dfrac{\partial \phi}{\partial \nu} =0,&\ \ {\rm on }\, \partial\Omega\, ,
\end{array}
\right.
\end{equation}
which  admits a divergent sequence $\lambda_n^{{\scriptscriptstyle\mathcal{N}}}$, $n\in \N$, of non-negative eigenvalues of finite multiplicity, with $\lambda_1^{{\scriptscriptstyle\mathcal{N}}}=0$.
The second is the eigenproblem
\begin{equation}\label{magbc}
\left\{
\begin{array}{ll}
\cu\, \cu \bm{\psi} =\lambda \bm{\psi} ,& \ \ {\rm in}\ \Omega ,\\
\di \bm{\psi}=0& \ \ {\rm in}\ \Omega ,\\
\bm{\nu} \times \cu \bm{\psi}  =\bm{0} ,& \ \ {\rm on}\ \Gamma ,\\
 \bm{\psi} \cdot \bm{\nu}  =0 ,& \ \ {\rm on}\ \Gamma
\end{array}
\right.
\end{equation}
which   admits a divergent sequence $\lambda_n^{{\scriptscriptstyle\mathcal{M}}}$, $n\in \N$, of non-negative eigenvalues of finite multiplicity.

 In the  following result from \cite[Thm.~3.10]{lamstr},  one has actually to better assume the condition $k\ne 0$. On the other hand  it is straightforward that  if $\Omega$ is simply connected and $k=0$ then $0\notin \Sigma$ because  the corresponding eigenfunctions would have zero  div and zero curl (see \eqref{minmax}) hence, being tangential, they would be identically zero, see   \cite[Prop.~2,~p.~219]{dali3} for more information; see also the more recent paper \cite{ciarlet}.

\begin{theorem}   Assume that $k\ne 0$ and $\theta >0$.   We have that $0\in \Sigma $ if and only if  $k^2 \in \{\theta \lambda_n^{{\scriptscriptstyle\mathcal{N}}}:\ n\in \N \} \cup  \{ \lambda_n^{{\scriptscriptstyle\mathcal{M}}}:\ n\in \N  \}$.
\end{theorem}

\subsection{Remarks on trace problems and Steklov expansions}

We denote by   $\bm{E}_n^{\scriptscriptstyle\Omega}$, $n\in \N$,  an orthonormal sequence of eigenvectors associated with the eigenvalues $\lambda_n$ of problem  \eqref{classiccucu1},  where it is understood that they are normalized with respect to the quadratic from $\mathcal{Q}$.

Let  $\pi_{\rm \scriptscriptstyle T}$ denote the trace operator from $X_{\rm \scriptscriptstyle T}(\Omega)$ to $TL^2(\Gamma)$ where
$$TL^2(\Gamma) =  \{{\bm u} \in (L^2(\Gamma))^3 : {\bm \nu} \cdot {\bm u} = 0 \,\, {\rm on}\ \Gamma \} . $$
By setting
 \begin{equation}\label{norm}
\bm{E}_n^{\scriptscriptstyle\Gamma}:= \sqrt{|\lambda_n-\eta |}\, \pi_{\rm \scriptscriptstyle T}\bm{E}_n^{\scriptscriptstyle\Omega}
\end{equation}
one can prove, in the spirit of Section~\ref{tracesec}, that $\bm{E}_n^{\scriptscriptstyle\Gamma}$, $n\in \N$, is an orthonormal basis of $TL^2(\Gamma)$.

These bases can be used to represent the solutions of the following problem
\begin{equation}\label{classiccucu10}
\left\{
\begin{array}{ll}
\cu\, \cu \bm{U} -k^2 \bm{U}  -\theta\, \gr\,  \di \bm{U} =\bm{0} ,& \ \ {\rm in}\ \Omega ,\\
\bm{\nu} \times \cu \bm{U}  = \bm{f}  ,& \ \ {\rm on}\ \Gamma ,
\end{array}
\right.
\end{equation}
where $ \bm{f} \in TL^2(\Gamma)$.

Let $ \bm{f}$ have the following representation
 \begin{equation*}
 \bm{f}=\sum_{n=1}^{\infty }c_n \bm{E}_n^{\scriptscriptstyle\Gamma},
\end{equation*}
with $(c_n)_{n\in \N}\in \ell^2$.
It is proved in \cite[Thm.~4.1]{lamstr} that if $0\notin \Sigma$ then the solution $ \bm{U}$ of \eqref{classiccucu10} can be expanded in terms of the above basis as follows
\begin{equation}\label{representationu}
\bm{U}= \sum_{n=1}^{\infty } \left(   \frac{\sqrt{|\lambda_n-\eta|}}{\lambda_n}\, c_n \right) \bm{E}_n^{\scriptscriptstyle\Omega}\, .
\end{equation}

Finally, under our assumptions we can represent the trace space of   $X_{\rm \scriptscriptstyle T}(\Omega)$ as follows

\begin{equation}\pi_{\rm \scriptscriptstyle T} (  X_{\rm \scriptscriptstyle T}(\Omega)     )=
 \pi_{\rm \scriptscriptstyle T} ( \mathcal{H}(\Omega)) =
\left\{  \sum_{j=1}^{\infty} c_j  \bm{E}_j^{\scriptscriptstyle\Gamma}: \ \sum_{j=1}^{\infty}|\lambda_j-\eta ||c_j|^2< \infty  \right\}
 \end{equation}
 which is the counterpart of the representation \eqref{auchmuty} for our problem.

\section{The case where $\Omega $ is the unit ball}
\label{ballsec}

In this section we consider problem \eqref{classiccucu1} in $\Omega = B$, where $B$ is the unit ball in $\R^3$ centred at zero and  we compute explicitly its eigenvalues and eigenvectors. In particular, we shall prove that there exist two families of eigenvectors, one of which is not divergence-free.

We proceed to define the vector spherical harmonics, following the notation of \cite{kri}. Recall that the (scalar) spherical harmonics are given by
\begin{equation} \label{sph_h}
Y_{\sigma m l}(\theta, \phi) = \sqrt{\frac{\eps_m}{2 \pi}} \sqrt{\frac{(2l+1)(l-m)!}{2 (l+m)!}}P^m_l(\cos \theta) F_{\sigma}(\phi) = C_{lm} P^m_l(\cos \theta ) F_{\sigma}(\phi)
\end{equation}
where $\theta \in [0, \pi]$, $\phi \in [0, 2 \pi)$, $\sigma \in \{\rm e, o\}$, $l \in \N \cup \{0\}$, $m \in \N \cup \{0\}$, $m \leq l$, $P^m_l$ is the Legendre function associated with the Legendre polynomial $P_l$, $\eps_m = 2 - \delta_{m0}$, and
\[
F_{\rm e}=\cos m \phi, \ \      F_{\rm o}= \sin m \phi .
\]
According to \cite[p.~627]{kri} we will use the multi-index notation
\[Y_n := Y_{\sigma m l},
 \]
 and we shall denote by  $\mathscr{O} $ the set of  triple indices $\sigma m l$  where $\sigma \in \{\rm e, o\}$,   $l\in \N \cup \{0\}$ and  $m\in \{0, \dots, {\it  l}  \}$.
For $n \in \mathscr{O}$, let
\begin{align}\label{def: vec sph h}
&\bs{A}_{1n}(\bs{\xi})= \frac{1}{\sqrt{l (l+1)}}\gr_{\bs{\xi}} Y_{n}(\bs{\xi}) \times \bs{\xi} \\
&\bs{A}_{2n}(\bs{\xi}) = \frac{1}{\sqrt{l (l+1)}} \gr_{\bs{\xi}} Y_{n}(\bs{\xi}) \\
&\bs{A}_{3n}(\bs{\xi}) = Y_{n}(\bs{\xi}) \bs{\xi}
\end{align}
for all $\bs{\xi} \in \p B$, where $Y_{n} = Y_{\sigma m l}$ is defined as in \eqref{sph_h}, see \cite[p.~350]{kri}. We extend this definition to points $\bm{x} \in B \setminus \{\bm{0}\}$ by setting
\begin{align}\label{def: vec sph h 2}
&\bs{A}_{1n}(\bm{x})=  \frac{|\bm{x}|}{\sqrt{l (l+1)}} \left( \gr_{\bm{x} }Y_{n}\bigg(\frac{\bm{x}}{|\bm{x}|}\bigg) \times \frac{\bm{x}}{|\bm{x}|}\right)
\\
&\bs{A}_{2n}(\bm{x}) = \frac{|\bm{x}|}{\sqrt{l (l+1)}} \gr_{\bm{x}} Y_{n}\bigg(\frac{\bm{x}}{|\bm{x}|}\bigg)  \\
&\bs{A}_{3n}(\bm{x}) = Y_{n}\bigg(\frac{\bm{x}}{|\bm{x}|}\bigg) \frac{\bm{x}}{|\bm{x}|}
\end{align}
By definition $\bs{A}_{1\sigma 0 0} = \bs{A}_{2\sigma 0 0} = 0$. The family $\{ \bs{A}_{\tau n}:\   \tau\in\{1,2,3 \},\   n \in \mathscr{O}\}$ is a complete orthonormal system in $L^2(\p B)^3$. Therefore, we can expand any vector field $\bs{E} \in L^2(B)^3$ as follows:
\[
\bs{E}(r, \bs{\xi}) = \sum_{n \in \mathscr{O}} \bigg( E^1_{n}(r) \bs{A}_{1n}(\bs{\xi}) + E^2_{n}(r) \bs{A}_{2n}(\bs{\xi}) + E^3_{n}(r) \bs{A}_{3n}(\bs{\xi}) \bigg).
\]
We note immediately that since we are interested in solutions of \eqref{classiccucu1}, the boundary condition $\bs{E} \cdot \bs{\nu} = 0$ is equivalent to $E^3_{n}(1) = 0$, since one easily checks that $\bs{A}_{2n}(\bs{\xi}) \cdot \bs{\xi} = 0$ and $\bs{A}_{1n}(\bs{\xi}) \cdot \bs{\xi} = 0$. \\
Recall that the vector Laplacian acts in the following way
\begin{equation} \label{formula 1}
\Delta (f(r) \bs{V}(\bs{\xi})) = \bigg(\frac{1}{r^2} \frac{\p}{\p r}  \bigg(  r^2 \frac{\p f(r)}{\p r} \bigg)     \bigg) \bs{V}(\bs{\xi}) + f(r) \Delta \bs V(\bs{\xi}).
\end{equation}
We have the following formulae for the vector Laplacian of the vector spherical harmonics $\bs{A}_n$
\begin{equation}\label{formula 2}
\begin{aligned}
&\Delta \bs{A}_{1n} = - \frac{1}{r^2} l (l+1) \bs{A}_{1n} \\
&\Delta \bs{A}_{2n} = \frac{2}{r^2} \sqrt{(l(l+1))} \bs{A}_{3n} - \frac{1}{r^2} l (l+1) \bs{A}_{2n} \\
&\Delta \bs{A}_{3n} = - \frac{1}{r^2} (2 + l(l+1)) \bs{A}_{3n} + \frac{2 \sqrt{l (l+1)}}{r^2} \bs{A}_{2n}.
\end{aligned}
\end{equation}
Moreover, one can compute (see also (C.14) in \cite{kri})
\begin{equation} \label{formula 3}
\di \bs E(r, \bs{\xi}) = \sum_{n \in \mathscr{O}} \bigg( \frac{\p E^3_{n}(r)}{\p r} + \frac{2}{r} E^3_{n}(r) - \frac{\sqrt{l(l+1)}}{r} E^2_{n}(r)\bigg) Y_{n}(\bs{\xi}).
\end{equation}
Let us define
\[
\Phi(r) := \bigg( \frac{\p E^3_{n}(r)}{\p r} + \frac{2}{r} E^3_{n}(r) - \frac{\sqrt{l(l+1)}}{r} E^2_{n}(r)\bigg).
\]
Note that $\di \bs{E} = 0$ is equivalent to $\Phi = 0$.

\begin{remark}
\label{rem:family_sol}
In \cite{calamo}, the authors find two families of solutions to the equation $\cu \cu u - k^2 u = 0$: one is given by $\bs{M}_n = \cu (\bs{x} j_l(k |\bm{x}|) Y_{lm}(\bm{x}/|\bm{x}|))$, the other one by $\bs{N}_n = \cu \bs{M}_n$.   Here $j_l$ is the spherical Bessel function of the first kind of order $l$, namely  $j_l(z)= \sqrt{\pi/(2z)}J_{l+1/2}(z)$.
These two families are divergence-free by definition. Indeed, the solutions $\bs{M}_n$ have $E^j = 0$, $j=2,3$, while the solutions $\bs{N}_n$ satisfy $\Phi(r) = 0$ for $E^2$, $E^3$ not identically zero. The function $\bs{N}_n = \cu \bs{M}_n(\bm{x})$ in \cite{calamo}, for $n \in \mathscr{O}$, are given by
\[
\begin{split}
\cu \bs{M}_n(\bm{x}) &= \cu \left[j_l(k |\bm{x}|)\left( \gr_{\bm{x}} Y_{n}\bigg(\frac{\bm{x}}{|\bm{x}|}\bigg) \times \frac{\bm{x}}{|\bm{x}|}\right)\right] \\
&= \frac{j_l(k |\bm{x}|) l (l+1)}{|\bm{x}|} \bs{A}_{3n}  + \sqrt{l(l+1)}\bigg(j_l'(k |\bm{x}|) k + j_l(k |\bm{x}|) \frac{1}{|\bm{x}|} \bigg) \bs{A}_{2n} .
\end{split}
\]
This is consistent with the fact that $\di \cu \bs{M}_n = 0$.
\end{remark}
\vspace{0.2cm}
From \eqref{formula 3} it is easy to compute
\begin{equation} \label{formula 4}
\gr(\di \bs E)(r, \bs{\xi}) = \sum_{n \in \mathscr{O}} \Phi'(r) \bs{A}_{3n}(\bs{\xi}) + \sum_{n \in \mathscr{O}}  \frac{\Phi(r)}{r} \sqrt{l (l+1)} \bs{A}_{2n}(\bs{\xi})  .
\end{equation}
Due to \eqref{formula 1}, \eqref{formula 2}, \eqref{formula 4} we then have that
{\small
\begin{multline}\label{trecomponenti}
- \Delta \bs{E} + (1- \theta) \gr \di \bs{E} = \\
\sum_{n \in \mathscr{O}}\begin{pmatrix}
- \bigg( \frac{1}{r^2} \frac{\p}{\p r} r^2 \frac{\p E^1_{n}}{\p r} - \frac{E^1_{n} l (l+1)}{r^2} \bigg) \\
\bigg(- E^3_{n}\frac{2 \sqrt{l(l+1)}}{r^2} - \frac{1}{r^2} \frac{\p}{\p r} r^2 \frac{\p E^2_{n}}{\p r} + \frac{E^2_{n}l (l+1)}{r^2} + (1- \theta) \sqrt{l(l+1)}\frac{\Phi}{r} \bigg) \\
\bigg( - \frac{1}{r^2} \frac{\p}{\p r} r^2 \frac{\p E^3_{n}}{\p r} + \frac{(2 + l(l+1))E^3_{n}}{r^2} - E^2_{n} \frac{2 \sqrt{l (l+1)}}{r^2} + (1-\theta) \Phi'  \bigg)
\end{pmatrix} \cdot
\begin{pmatrix}
\bs{A}_{1n} \\
\bs{A}_{2n} \\
\bs{A}_{3n}
\end{pmatrix}
\end{multline}}
Consider now $- \Delta \bs{E} + (1-\theta) \gr \di \bs{E} - k^2 \bs{E} = 0$ in $B$. It is clear that the first equation   provided by the first row in \eqref{trecomponenti} above can be solved independently of the other two, by means of separation of variables and classical Sturm-Liouville theory, yielding
\begin{equation} \label{eq:E1}
E^1_{n}(r) = E^1_l(r) = j_l(k r)
\end{equation}
where $j$ is the spherical Bessel function of the first kind of order $l$ as above. The other two equations form a coupled system of Sturm-Liouville equations. For the sake of clarity, we first solve the system in the case $\theta = 1$.\\[0.3cm]
\textbf{Case $\theta = 1$}. We need to solve the two-parameter family of ODE systems
\begin{equation}
\label{sysSL}
\begin{cases}
    - \frac{1}{r^2} \frac{\p}{\p r} r^2 \frac{\p E^2_{n}}{\p r} +   \frac{   l (l+1) E^2_{n}  }{r^2}   -     \frac{2 \sqrt{l(l+1)}  E^3_{n}   }{r^2}           -k^2E^2_n = 0, \quad &\textup{in $(0,1),$} \vspace{2mm}\\
- \frac{1}{r^2} \frac{\p}{\p r} r^2 \frac{\p E^3_{n}}{\p r} + \frac{(2 + l(l+1))E^3_{n}}{r^2} -  \frac{2 \sqrt{l (l+1)}  E^2_{n}   }{r^2} -k^2E^3_n= 0, \quad &\textup{in $(0,1)$,} \vspace{2mm}\\
E^3_{n}(1) = 0. &\textup{}
\end{cases}
\end{equation}
Recall the spherical Bessel equation
\begin{equation}
\label{eq: sphBeq}
\frac{\p}{\p r}\left( r^2 \frac{\p}{\p r} f \right) + (k^2 r^2 - l (l+1)) f = 0.
\end{equation}
Define the spherical Bessel operator of indices $k \in \R$ and $l \in \Z$ as follows
\[
\cL_{k, l}(f) = \frac{\p}{\p r}\left( r^2 \frac{\p}{\p r} f \right) + (k^2 r^2 - l (l+1)) f.
\]
We will require $f \in H^2([0,1])$, which is equivalent to imposing that the solution $f$ to \eqref{eq: sphBeq} is not singular in zero. System \eqref{sysSL} can be then rewritten in the simpler form
\begin{equation}
\label{sysSL2}
\begin{cases}
&\cL_{k, l}(E^2) = - 2 \sqrt{l (l+1)} E^3, \vspace{2mm}\\
&\cL_{k, l}(E^3) = 2 E^3 - 2 \sqrt{l (l+1)} E^2, \vspace{2mm}\\
&E^3(1) = 0.
\end{cases}
\end{equation}
where we have omitted the dependence on $n$ in $E^2$, $E^3$.

\begin{remark}
Note that, upon replacing $E^2$ in the first equation, by means of the equality $ \frac{2E^3-\cL_{k, l}(E^3)}{ 2 \sqrt{l (l+1)}} = E^2$, we obtain the fourth-order equation
\[
\left(\cL_{k, l}\right)^2 (E^3) - 2 \cL_{k, l}(E^3) = 4 l (l+1) E^3,
\]
with the boundary condition $E^3(1) = 0$.
\end{remark}

In view of Remark \ref{rem:family_sol} a solution to the differential equation in \eqref{sysSL2} is given by
\begin{equation} \label{Monksol}
E^2(r) = \sqrt{l(l+1)} \bigg( j_l'(k r) k + \frac{j_l(k r)}{r} \bigg), \quad E^3(r) = \frac{j_l(k r) l (l+1)}{r};
\end{equation}
note however that the condition $E^3(1) = 0$ is not satisfied for general $k$. To overcome this problem, it is convenient to first find another explicit solution of \eqref{sysSL2}.

We claim that the functions
\begin{equation}
\label{def:cE}
\cE^2(r) := \sqrt{l(l+1)}\frac{j_l(kr)}{r}, \quad \cE^3(r) := k j'_l(k r),
\end{equation}
are solutions of the system of differential equations \eqref{sysSL2}. In the sequel we use the abbreviated  notation  $\cL := \cL_{k, l}$. Let us compute
\begin{equation} \label{conti1}
\cL(\!k\, j'_l(k\, r)) = k^{3}\, r^2\, j_l^{(3)}(k\, r) + 2k^2 \,r\, j_l''(k\, r) + k\,(k^2 r^2 - l (l+1))\, j_l'(k\, r)
\end{equation}
Recall that $j_l(\cdot)$ satisfies \eqref{eq: sphBeq}. Differentiation of \eqref{eq: sphBeq} with $f(r) = j_l(k r)$ gives
\begin{equation}\label{conti2}
k^{3}\,r^2\, j_l^{(3)}(k\, r) + 4 k^2 r\, j_l''(k\, r) + k( k^2 r^2 + 2 - l(l+1))\,j_l'(k r) + 2 k^2 r \, j_l(k \, r) = 0
\end{equation}
Equation \eqref{conti2} implies that we can rewrite \eqref{conti1} as follows
{\small
\begin{equation} \label{conti3}
\begin{split}
\cL(k j'_l(k r)) &= -4k^2 r j_l''(k r) - k( k^2 r^2 + 2 - l(l+1))j_l'(k r) \\ &\quad\quad\quad -2 k^2 r j_l(k r) + 2k^2 r j_l''(k r) + k (k^2 r^2 - l (l+1)) j_l'(k r) \\
&= -2 k^2 r j_l''(k r) - 2 k j_l'(k r) - 2 k^2 r j_l(k r)
\end{split}
\end{equation}
\small}
The spherical Bessel equation $\cL(j_l) = 0$, see \eqref{eq: sphBeq}, again implies that
\[
k^2 r \, j''_l(k r) = - 2 k j_l'(k r) - \big(k^2 r^2 - l(l+1)\big) \frac{j_l(k r)}{r}
\]
hence \eqref{conti3} can be rewritten as
{\small
\[
\begin{split}
\cL\left( k\, j'_l(k r)\right) &= 2 \left( 2 k j_l'(k r) + \big(k^2 r^2 - l(l+1)\big) \frac{j_l(k r)}{r}\right) - 2 k j_l'(k r) - 2 k^2 r \, j_l(k r) \\
&=2 k j_l'(k r) - \frac{2l (l+1)}{r} j_l(k r) = 2 \cE^3 - 2\sqrt{l (l+1)} \cE^2,
\end{split}
\]}
\!which implies the identity $\cL (\cE^3) = 2 \cE^3 - 2\sqrt{l (l+1)} \cE^2$. \\
We note that the functions $E^2, E^3$ defined in \eqref{Monksol}
satisfy the equalities
$$\cE^2 = E^3/\sqrt{l(l+1)},\  \ {\rm and} \ \ \cL(E^3) = 2E^3 - 2 \sqrt{l (l+1)} E^2.$$
 Then  it is immediate to check that
$
\cL(\cE^2) = - 2 \sqrt{l(l+1)}\cE^3
$
which implies  that the couple $(\cE^2, \cE^3)^t$ solves the two equations in \eqref{sysSL2}, as claimed.

Note that the divergence of the function
\[
\cE(r, \bs{\xi}) = \sum_{n \in \mathscr{O}} \bigg(\cE^2_{n}(r) \bs{A}_{2n}(\bs{\xi}) + \cE^3_{n}(r) \bs{A}_{3n}(\bs{\xi})\bigg)
\]
is non-trivial unless $k = 0$. Indeed,
\begin{equation} \label{eq:DivEcomp}
\begin{split}
\di \bs \cE(r, \bs{\xi}) &= \sum_{n \in \sO} \bigg( \frac{\p \cE^3(r)}{\p r} + \frac{2}{r} \cE^3(r) - \frac{\sqrt{l(l+1)}}{r} \cE^2(r)\bigg) Y_{n}(\bs{\xi})\\
&= \sum_{n \in \sO} \bigg( j_l''(k r) k^2 + \frac{2}{r} j_l'(k r) k - \frac{l(l+1)}{r} \frac{j_l(k r)}{r}\bigg) Y_{n}(\bs{\xi})\\
&= -k^2 \sum_{n \in \sO} j_l(k r) Y_{n}(\bs{\xi}) \neq 0,
\end{split}
\end{equation}
where in the last equality we used again \eqref{eq: sphBeq}.\\
We note {\it en passant} that the previous formula for the divergence, namely
\[
\bigg( \frac{\p \cE^3(r)}{\p r} + \frac{2}{r} \cE^3(r) - \frac{\sqrt{l(l+1)}}{r} \cE^2(r)\bigg) = - k^2 j_l(k r)
\]
gives the equality
\begin{equation} \label{ansatz_theta1}
\cE^2(r) = \frac{r}{\sqrt{l(l+1)}}\left( \frac{1}{r^2}\frac{\p (r^2 \cE^3(r))}{\p r} + k^2 j_l(k r) \right).
\end{equation}
A similar computation for the functions $E^2$, $E^3$ defined in \eqref{Monksol} gives
\begin{equation}
\label{ansatz_theta1_2}
E^2(r) = \frac{1}{\sqrt{l(l+1)} \, r} \frac{\p (r^2 E^3(r))}{\p r}
\end{equation}
in agreement with Formula (7.4) in \cite{kri}. \\
For the boundary condition to be satisfied we then choose a linear combination of $(E^2, E^3)^t$ and $(\cE^2, \cE^3)^t$. If $j'_l(k) = 0$ then $(\cE^2, \cE^3)^t$ is already a solution of \eqref{sysSL2} with non-trivial divergence. Otherwise, we just set
\begin{equation}
\label{eq:F}
\begin{aligned}
F^2(r) &:= - \frac{(l(l+1))^{3/2}j_l(k)}{k j_l'(k)} \frac{j_l(k r)}{r} + \sqrt{l(l+1)}\left( \frac{j_l(k r)}{r} + k j_l'(k r) \right) \\
F^3(r) &:= - l(l+1) \bigg(\frac{j_l(k)}{j_l'(k)} j_l'(kr) -  \frac{j_l(k r)}{r} \bigg),
\end{aligned}
\end{equation}
and then $\bs{F}:=(F^2, F^3)^t$ solves \eqref{sysSL2} with the right boundary condition. For subsequent  use, note that
\[
\begin{split}
F^2(1) &= \sqrt{l(l+1)}\bigg(1 - \frac{l(l+1)j_l(k)}{k j_l'(k)} \bigg)j_l(k)+ \sqrt{l(l+1)}k j_l'(k) \\
&= \sqrt{l(l+1)} \frac{j_l(k) j_l'(k) k - l(l+1) j_l(k)^2 + k^2 (j_l'(k))^2}{k j_l'(k)}\\
\end{split}
\]
for all $l \geq 1$.\\
Let us also compute
{\small
\[
\begin{split}
&\frac{1}{\sqrt{l(l+1)}}\frac{\d}{\d r} F^2(r) = \bigg(k \frac{j_l'(kr)}{r} - \frac{j_l(k r)}{r^2}\bigg) \bigg( 1 - \frac{l (l+1) j_l(k)}{k j_l'(k)} \bigg) + k^2 j_l''(k r) \\
&\overset{\eqref{eq: sphBeq}}{=}  \frac{j_l(k r)}{r^2} \bigg( \frac{l (l+1) j_l(k)}{k j_l'(k)} - k^2 + \frac{l (l+1) - 1}{r^2} \bigg) - k \frac{j_l'(kr)}{r} \bigg( 1 + \frac{l (l+1) j_l(k)}{k j_l'(k)} \bigg)
\end{split}
\]
}
therefore
\[
\bigg(\frac{\d}{\d r} F^1(r) \bigg)|_{r = 1} = \sqrt{l(l+1)} \bigg( - (k^2 + 1) j_l(k) - k j_l'(k) + \frac{j_l(k)^2 l (l+1)}{k j_l'(k)} \bigg)
\]
Note that the solution $\bs{E}$ of
\[
\begin{cases}
\cu \cu \bs{E} - k^2 \bs{E} = 0, &\textup{in $B$,} \\
\bs{\nu} \times \cu \bs{E} - \la (\bs{\nu} \times \bs{E} \times \bs{\nu}) = 0, &\textup{on $\p B$,} \\
\bs{\nu} \cdot \bs{E} = 0, &\textup{on $\p B$,}
\end{cases}
\]
can be found in the form
\[
\bs{E} =  \sum_{n \in \sO} a_n (F_n^2 \bs{A}_{2n} + F_n^3 \bs{A}_{3n}) + b_n E^1_{l} \bs{A}_{1n}
\]
where $F^2_n, F^3_n$ are defined in \eqref{eq:F}, $E^1_n$ is defined in \eqref{eq:E1}, and $a_n, b_n \in \C$.
For the following computations it is useful to recall the formulae (see also  \cite[p.~633]{kri})
\begin{align*}
&\cu(f(r) \bs{A}_{1n}) = \bigg(\frac{\sqrt{l(l+1)}}{r} f(r) \bigg) \bs{A}_{3n} + \bigg(\frac{1}{r} \frac{\p}{\p r}\left( r f(r) \right)\bigg) \bs{A}_{2n}, \\
&\cu(f(r) \bs{A}_{2n}) = - \bigg(\frac{1}{r} \frac{\p}{\p r} (r f(r) )\bigg) \bs{A}_{1n}, \\
&\cu(f(r)\bs{A}_{3n}) = \bigg( \frac{\sqrt{l(l+1)}}{r} f(r) \bigg) \bs{A}_{1n}.
\end{align*}
In particular, we have that
\[
\begin{split}
&(\bs{\xi} \times \cu (F_n^2 \bs{A}_{2n} + F_n^3 \bs{A}_{3n}))|_{r = 1} = \\
&- \bs{\xi} \times ((F_n^2)'(1) + F_n^2(1))\bs{A}_{1n} + \bs{\xi} \times (\sqrt{l(l+1)} F^3_n(1) \bs{A}_{1n})\\
&= - ((F_n^2)'(1) + F_n^2(1))\bs{A}_{2n}
\end{split}
\]
Imposing the Steklov condition and taking into account that $\bs{\nu} \cdot \bs{E} = 0$ give
\[
\begin{split}
0 = (\bs{\nu} \times \cu \bs{E} &- \la (\bs{\nu} \times \bs{E} \times \bs{\nu}))|_{r = 1}\\
&= \sum_{n \in \sO} - a_n ((F_n^2)'(1) + (1 + \la) F_n^2(1)) \bs{A}_{2n} \\
& \hspace{2cm} \sum_{n \in \sO} - b_n \big[j_l(k) + j_l'(k) k + \la j_l(k)\big] \bs{A}_{1n}.
\end{split}
\]
Hence, the eigenvalues $\la_n$ of \eqref{classiccucu1} are given by two families. The first one is obtained by imposing
\[
\begin{split}
&0 =((F_n^2)'(1) + (1 + \la) F_n^2(1))\\
&\Rightarrow \la^{(1)}_l = \frac{ k^2 \,k\, j_l(k)\, j_l'(k)}{k j_l(k) j_l'(k) - l (l+1) (j_l(k))^2 + k^2 (j_l'(k))^2}
\end{split}
\]
for all $l \geq 1$. The second one is obtained by imposing that the coefficient of $\bs{A}_{1n}$ vanishes for some $n \in \sO$. We then obtain
\begin{equation} \label{eigen2}
\la^{(2)}_l = - \frac{j_l(k) + j_l'(k) k}{j_l(k)}
\end{equation}
It is not difficult to check that the eigenvalues $\la^{(2)}_l \distas{l \to +\infty} -l$.

Regarding $\la^{(1)}_l$, we note that, due to the recurrence formulae for the derivatives of Bessel functions,
\[
\frac{k j_l(k) j_l'(k) - l (l+1) (j_l(k))^2 + k^2 (j_l'(k))^2}{k^2} = - j_{l-1}(k)j_{l+1}(k)
\]
hence,
\begin{equation}
\label{la1_new}
\la^{(1)}_l = - \frac{k j_l(k) j_l'(k)}{j_{l+1}(k)j_{l-1}(k)}
\end{equation}
for $l \geq 1$. Now, recalling that $j_l(z) = \sqrt{\frac{\pi}{2z}} J_{l + 1/2}(z)$, where $J$ is the Bessel function of the first kind, and that we have the following large index asymptotic formula:
\[
J_{s}(z) \distas{s \to +\infty} \frac{1}{\sqrt{2 \pi s}} \left(\frac{\e z}{2 s}\right)^{s}
\]
for all $z \in \C$, we deduce that
\begin{equation}
\label{asymp_j_l}
j_l(z) \distas{l \to \infty} \frac{\e^{l + 1/2} z^l}{2^{l+3/2} (l + 1/2)^{l+1}}
\end{equation}
for $l \geq 1$. By using formulae \eqref{la1_new} and \eqref{asymp_j_l} we obtain
\[
\la^{(1)}_l \distas{l \to \infty} - l
\]
therefore the eigenvalues of the first family are diverging to $-\infty$ as $l \to +\infty$ (as expected).\\

\begin{figure}
\centering
\includegraphics[width=0.8\linewidth]{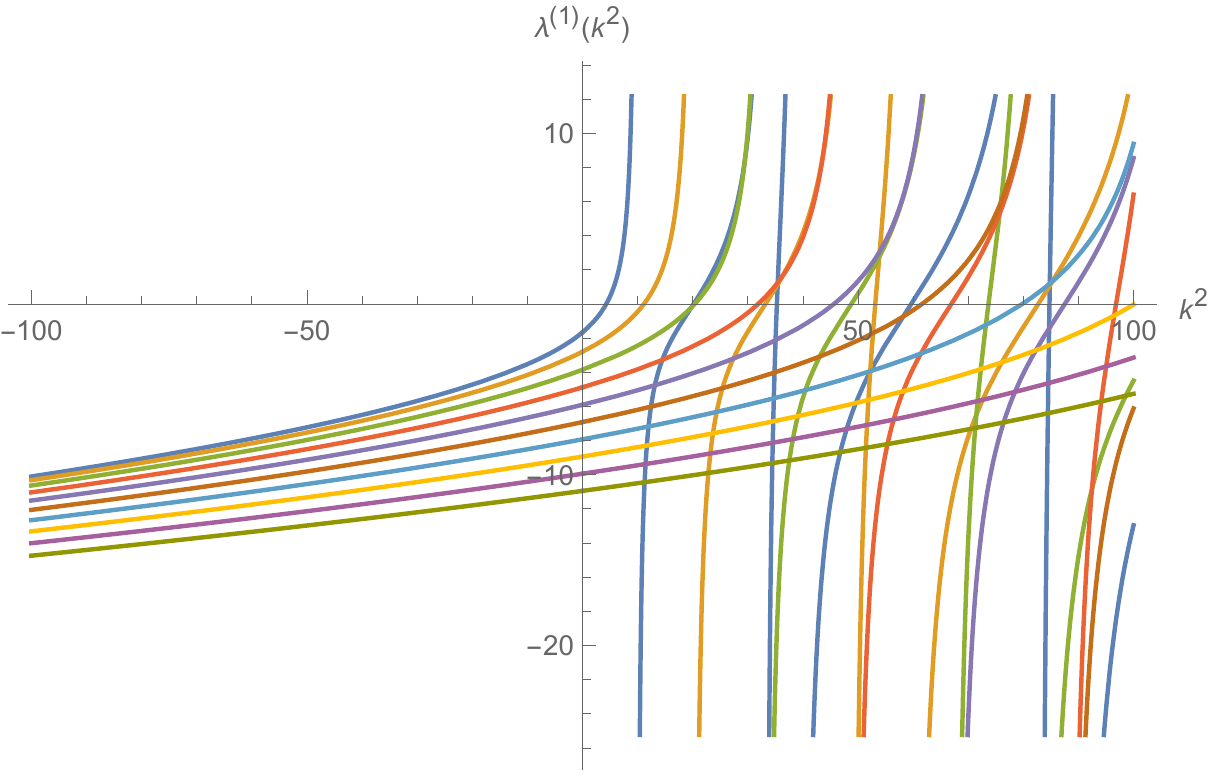}
\caption{Eigenvalues $\la^{(1)}_l$, $l=1, \dots, 10$, as a function of $k^2 \in (-100,100)$.}
\end{figure}

\noindent \textbf{General case $0 < \theta < + \infty$}. Since we are interested in solutions of
\begin{equation} \label{maineq}
\cu \cu \bs{E} - \theta \gr \di \bs{E} - k^2 \bs{E} = 0
\end{equation}
with non-trivial divergence, it is convenient to find an equation for $\di \bs{E}$.  The following discussion applies to a general domain $\Omega$ (not necessarily a ball).  Let $\varphi \in H^2(\Omega)$ be an arbitrary function.  Multiply \eqref{maineq} by $\gr \varphi$ and integrate in $\Omega$ to obtain
\[
\int_{\Omega} (\cu \cu \bs{E} \cdot \gr \varphi - \theta\, \gr \di \bs{E} \cdot \gr \varphi - k^2 \, \bs{E} \cdot \gr \varphi )\,dx= 0.
\]
Due to the Steklov boundary condition, the identity $\cu \gr \varphi = 0$, and standard integration by parts we have
\[
\begin{split}
\int_{\Omega} \cu \cu \bs{E} \cdot \gr \varphi \,dx= \int_{\Gamma} (\bs{\nu} \times \cu \bs{E}) \cdot \gr_{\Gamma} \varphi \, d\sigma  \\
= \la \int_{\Gamma} \bs{E} \cdot \gr \varphi \, d\sigma =  \la \int_{\Gamma} (\di_{\Gamma} \bs{E}) \, \varphi \, d\sigma ,
\end{split}
\]
where $\gr_{\Gamma} $ and $\di_{\Gamma}$ denote the standard tangential gradient and tangential divergence respectively.  Moreover,
\[
- \theta \int_{\Omega} \gr \di \bs{E} \cdot \gr \varphi\, dx= \theta \int_{\Omega} \Delta \di \bs{E} \varphi \,dx - \theta \int_{\Gamma} \frac{\p}{\p \nu}(\di \bs{E}) \varphi \, d\sigma ,
\]
as well as
\[
- k^2 \int_{\Omega} \bs{E} \cdot \gr \varphi  \, dx= k^2 \int_{\Omega} \di \bs{E} \varphi \,  dx - k^2 \int_{\Gamma} \underbrace{(\bs{E} \cdot \bs{\nu})}_\text{$= 0$} \varphi\, d\sigma ,
\]
hence
\[
\int_{\Omega} (\theta \Delta \di \bs{E} + k^2 \di \bs{E}) \varphi\, dx= \int_{\Gamma} \left( \la \di_{\Gamma} \bs{E} +\theta \frac{\p}{\p \nu}     (\di \bs{E}) \right) \varphi \, d\sigma ,
\]
which    implies that $\di \bs{E}$   satisfies the following boundary value problem
\begin{equation}
\label{eq: sysDivE}
\begin{cases}
-\Delta \di \bs{E} - \frac{k^2}{\theta} \di \bs{E} = 0, \quad &\textup{in $\Omega$,} \vspace{2mm}\\
\frac{\p}{\p \bs\nu} \di \bs{E} = - \frac{\la}{\theta} \di_{\Gamma} \bs{E}, \quad &\textup{on $\Gamma$.}
\end{cases}
\end{equation}

Returning to the case of the ball,
the equation in $\Omega = B$ has the explicit solutions
\[\di \bs{E}_n ( r, \bs{\xi}) = a j_l\left( \frac{k}{\sqrt{\theta}} r \right) Y_{n}(\bs{\xi}), \quad n \in \sO,\  a\in \C  . \]
Therefore, arguing as in \eqref{eq:DivEcomp} we have the equality
\[
\left(\frac{\p \cE^3(r)}{\p r} + \frac{2}{r} \cE^3(r) - \frac{\sqrt{l(l+1)}}{r} \cE^2(r) \right) = a j_l\left( \frac{k}{\sqrt{\theta}} r \right)
\]
which leads to the {\it ansatz} (see also \eqref{ansatz_theta1})
\begin{equation}\label{eq:ansatzgen}
\cE^2(r) = \frac{r}{\sqrt{l(l+1)}}\left( \frac{1}{r^2}\frac{\p (r^2 \cE^3(r))}{\p r} - a j_l\left(\frac{k}{\sqrt{\theta}} r\right) \right).
\end{equation}
According to the previous computations and keeping in mind the  case $\theta = 1$, in view of \eqref{def:cE}  we claim that
\[
\cE^2(r) = \sqrt{l(l+1)}\, \frac{j_l\left( \frac{k}{\sqrt{\theta}} r \right)}{r}, \quad\quad \cE^3(r) = \frac{k}{\sqrt{\theta}}j_l'\left( \frac{k}{\sqrt{\theta}} r \right) ,
\]
verify \eqref{eq:ansatzgen} for $a = - k^2/ \theta$ and they are solutions of the system \eqref{sysSL}. Let us show that \eqref{eq:ansatzgen} is satisfied. Note that
\[
\begin{split}
\frac{1}{r^2}\frac{\p (r^2 \cE^3(r))}{\p r} &= \frac{\p \cE^3(r)}{\p r} + \frac{2}{r} \cE^3(r) \\
&= \frac{k^2}{\theta} j_l''\left( \frac{k}{\sqrt{\theta}} r \right) + \frac{2}{r} \frac{k}{\sqrt{\theta}} j_l'\left( \frac{k}{\sqrt{\theta}} r \right) \\
&=  \left(-\frac{k^2}{\theta}  + \frac{l (l+1)}{r^2} \right) j_l\left( \frac{k}{\sqrt{\theta}} r \right)
\end{split}
\]
where the last equality is deduced by \eqref{eq: sphBeq} with $k$ replaced by $k/\sqrt{\theta}$,
so the right-hand side of \eqref{eq:ansatzgen} becomes
\[
\frac{r}{\sqrt{l(l+1)}}\left(\bigg(-\frac{k^2}{\theta} - a\bigg)  + \frac{l (l+1)}{r^2}  \right) j_l\left( \frac{k}{\sqrt{\theta}} r \right)
\]
which equals $\cE^2(r) = \sqrt{l(l+1)} \frac{j_l\left( \frac{k}{\sqrt{\theta}} r \right)}{r}$ for $a = - k^2/ \theta$, as claimed. With similar computations one can check that  the couple $(\cE^2, \cE^3)$
satisfies the system
{\small
\begin{equation}
\label{sysSLgen}
\begin{cases}
\begin{aligned}
&- \frac{1}{r^2} \frac{\p}{\p r} r^2 \frac{\p E^2_{n}}{\p r} + \frac{l (l+1) E^2_{n}}{r^2} -  \frac{2 \sqrt{l(l+1)}  E^3_{n}     }{r^2} - k^2 E^2_{n}\\
&\hspace{0.7cm}+ \frac{(1- \theta)\, \sqrt{l (l+1)}}{r} \bigg( \frac{\p E^3_{n}(r)}{\p r} + \frac{2}{r} E^3_{n}(r) - \frac{\sqrt{l(l+1)} E^2_{n}(r)  }{r} \bigg)  = 0,
\end{aligned} \: &\textup{in $(0,1),$}\vspace{2mm}\\
\begin{aligned}
&- \frac{1}{r^2} \frac{\p}{\p r} r^2 \frac{\p E^3_{n}}{\p r} + \frac{(2 + l(l+1))E^3_{n}}{r^2} - \frac{2 \sqrt{l (l+1)}  E^2_{n}    }{r^2} - k^2 E^{3}_{n} \\
&\hspace{0.7cm} + (1- \theta) \frac{\p}{ \p r}\bigg( \frac{\p E^3_{n}(r)}{\p r} + \frac{2}{r} E^3_{n}(r) - \frac{\sqrt{l(l+1)} E^2_{n}(r) }{r} \bigg) = 0,
\end{aligned} \: &\textup{in $(0,1)$,}
\end{cases}
\end{equation}}
\noindent which replaces system \eqref{sysSL} used in the case $\theta =1$.  Clearly, another solution
(which is divergence-free) is given by  $(E_{n}^2, E_{n}^3)$ as defined in \eqref{Monksol}.

As in the case $\theta = 1$ we consider the linear combination $\binom{F_{n}^2}{F_{n}^3} := a \binom{E_{n}^2}{E_{n}^3} + b \binom{\cE_{n}^2}{\cE_{n}^3}$ and we impose  the boundary condition $F_{n}^3(1) = 0$, corresponding to $\bs{E} \cdot \bs{\nu} = 0$ on $\p B$.
We obtain
\[
a = - b \frac{j_l'\left( \frac{k}{\sqrt{\theta}}\right) \frac{k}{\sqrt{\theta}}}{j_l(k) l (l+1)},
\]
hence, up to a constant factor,
\begin{equation}
\label{F2}
F_{n}^2(r) := - \frac{j_l'\left( \frac{k}{\sqrt{\theta}}\right) \frac{k}{\sqrt{\theta}}}{j_l(k) \sqrt{l (l+1)}} \bigg(\frac{j_l(k r)}{r}  + j_l'(k r) k \bigg) + \sqrt{l(l+1)} \frac{j_l\left( \frac{k}{\sqrt{\theta}} r \right)}{r}
\end{equation}
\begin{equation}
\label{F3}
F_{n}^3(r) := - \frac{j_l'\left( \frac{k}{\sqrt{\theta}}\right) \frac{k}{\sqrt{\theta}}}{j_l(k)} \frac{j_l(k r)}{r} + j_l'\left( \frac{k}{\sqrt{\theta}} r \right) \frac{k}{\sqrt{\theta}}
\end{equation}
As in the case $\theta = 1$, we deduce that there are two families of eigenvalues diverging to $- \infty$: the first one coincides with $\la_n^{(2)}$ defined in \eqref{eigen2}, and it is associated with eigenfunctions in the form $\bs{E}_n(r, \bs{\xi}) = j_l(kr) \bs{A}_{1n}(\bs{\xi})$; the second one is obtained as in the case $\theta = 1$ by imposing the Steklov boundary conditions, and they are given explicitly as solutions of the equation
\[
\la_l^{(1)} = - \frac{F_l^2(1) + (F_l^2)'(1)}{F_l^2(1)}.
\]
We proceed therefore to computing $F_l^2(1)$ and $(F_l^2)'(1)$. To simplify the notation we will not write the index $l$. We easily obtain
\[
F^2(1) = \frac{j\!\left(\frac{k}{\sqrt{\theta}}\right) j(k)\, l (l+1) - j'\!\left( \frac{k}{\sqrt{\theta}}\right) j'(k)\, \frac{k^2}{\sqrt{\theta}} - j'\!\left( \frac{k}{\sqrt{\theta}}\right) j(k)\,\frac{k}{\sqrt{\theta}} }{j(k) \, \sqrt{l (l+1)}}
\]
Next, we compute
\[
\begin{split}
&\frac{\p}{\p r}F^2(r)
= - \frac{j'\left(\frac{k}{\sqrt{\theta}}\right)\frac{k}{\sqrt{\theta}}}{j(k) \sqrt{l (l+1)}}  \left( \frac{k j'(k r)}{r} - \frac{j(k r)}{r^2}  + j''(k r) k^2 \right) \\
&\hspace{5cm} + \sqrt{l(l+1)} \bigg( \frac{k}{\sqrt{\theta}} \frac{j'\left( \frac{k}{\sqrt{\theta}} r\right)}{r} - \frac{j\left( \frac{k}{\sqrt{\theta}}r\right)}{r^2} \bigg) \\
&\overset{\eqref{eq: sphBeq}}{=} - \frac{j'\left( \frac{k}{\sqrt{\theta}}\right)\frac{k}{\sqrt{\theta}}}{j(k) \sqrt{l (l+1)}}  \left( -\frac{k j'(k r)}{r} + \frac{j(k r)}{r^2} (l (l+1) - k^2 r^2 - 1) \right) \\
&\hspace{1cm} + \frac{1}{j(k) \sqrt{l (l+1)}} \bigg( \frac{k}{\sqrt{\theta}} l(l+1) j(k) \frac{j'\left( \frac{k}{\sqrt{\theta}} r\right)}{r} \bigg) -  \sqrt{l(l+1)}\frac{j\left( \frac{k}{\sqrt{\theta}}r\right)}{r^2} \\
\end{split}
\]
hence
\[
(F^2)'(1) = \frac{j'\!\left(\frac{k}{\sqrt{\theta}}\right) \frac{k}{\sqrt{\theta}} \big(k j'(k) + j(k) (k^2 + 1 )\big) - j\left( \frac{k}{\sqrt{\theta}}\right) j(k) l (l+1)}{j(k) \sqrt{l (l+1)}}
\]
Therefore,
\[
\begin{split}
\la_l^{(1)} &= -\frac{F_l^2(1) + (F_l^2)'(1)}{F_l^2(1)} \\
&= - \frac{j'\!\left(\frac{k}{\sqrt{\theta}}\right) j(k) \frac{k}{\sqrt{\theta}} \, k^2}{j\!\left(\frac{k}{\sqrt{\theta}}\right) j(k)\, l (l+1) - j'\!\left( \frac{k}{\sqrt{\theta}}\right) j'(k)\, \frac{k^2}{\sqrt{\theta}} - j'\!\left( \frac{k}{\sqrt{\theta}}\right) j(k)\,\frac{k}{\sqrt{\theta}} }
\end{split}
\]
which is in agreement with the case $\theta = 1$.

We summarize the previous discussion in the following  theorem where it is understood that the values of $k$ are such that the denominators in \eqref{labda1} and \eqref{labda2} do not vanish (in the remark below we explain the meaning of this condition).   

\begin{theorem} \label{eigen_final}
If $k\ne 0$ then the eigenvalues and eigenfunctions of the Steklov problem \eqref{classiccucu1} in the unit ball $B$ of $\R^3$ are given by the following two families:
\begin{equation}
\label{labda1}
\begin{cases}
\la_n^{(1)} = - \frac{j'_l\!\left(\frac{k}{\sqrt{\theta}}\right) j_l(k) \frac{k}{\sqrt{\theta}} \, k^2}{j_l\!\left(\frac{k}{\sqrt{\theta}}\right) j_l(k)\, l (l+1) - j_l'\!\left( \frac{k}{\sqrt{\theta}}\right) j_l'(k)\, \frac{k^2}{\sqrt{\theta}} - j_l'\!\left( \frac{k}{\sqrt{\theta}}\right) j_l(k)\,\frac{k}{\sqrt{\theta}} }, \quad &n \in \sO \vspace{2mm}\\
\bs{F}_n = F^2_n \bs{A}_{2n} +  F^3_n \bs{A}_{3n},\quad &\textup{$n \in \sO$,} \\
\end{cases}
\end{equation}
and
\begin{equation}
\label{labda2}
\begin{cases}
\la^{(2)}_n = - \frac{j_l(k) + j_l'(k) k}{j_l(k)}, \quad &n \in \sO \vspace{2mm}\\
\bs{E}_n = E^1_n \bs{A}_{1n}, \quad &n \in \sO \vspace{2mm}\\
\di \bs{E}_n = 0. &
\end{cases}
\end{equation}
where $\bs{A}_{\tau n}$, $\tau=1,2,3$, $n \in \sO$ are the vector spherical harmonics defined in \eqref{def: vec sph h}, the functions $F^2_n$, $F^3_n$ are defined respectively in \eqref{F2}, \eqref{F3} and $E^1_n$ is defined in \eqref{eq:E1}.

\end{theorem}

\begin{remark} The squared values of $k\ne 0$ for which the denominators in \eqref{labda1} and \eqref{labda2} vanish are the eigenvalues $A$ of the Dirichlet problem \eqref{dirichlet}. Indeed, our computations in this section show that the solutions of the partial differential equation in \eqref{classiccucu1}, which are tangential at the boundary of $B$, are given by the two families of vector fields
$\bs{F}_n = F^2_n \bs{A}_{2n} +  F^3_n \bs{A}_{3n}$ and  $\bs{E}_n = E^1_n \bs{A}_{1n}$. Note that $\bs{E}_n$ is tangential because $ \bs{A}_{1n}$ is tangential, while $\bs{F}_n$ is tangential because  $F^3_n(1)=0$.

Consider now the condition
\begin{equation}\label{secondden}
 E^1_n(1)=0\,.
\end{equation}
The values of $k$ for which \eqref{secondden} is satisfied are exactly the values of $k$ for which the denominator in \eqref{labda2} vanishes: in this case the vector field $\bs{E}_n$ vanishes at the boundary of $B$, hence $\bs{E}_n$ is a solution of \eqref{dirichlet} with $A=k^2$.

Assume now that the denominator in \eqref{labda2} does not vanish. In this case, the function $F^2_n$ is well-defined, hence the solution $\bs{F}_n$ is well-defined as well.  Under this assumption consider the condition
\begin{equation}\label{firstden}
 F^2_n(1)=0\,.
\end{equation}
The values of $k$ for which \eqref{firstden} is satisfied are exactly the values of $k$ for which the denominator in \eqref{labda1} vanishes: in this case the vector field $\bs{F}_n$ vanishes at the boundary of $B$, hence $\bs{F}_n$ is a solution of \eqref{dirichlet} with $A=k^2$.

 Thus, assuming that the denominators of \eqref{labda1} and   \eqref{labda2} do not vanish corresponds to our assumption \eqref{exclusion}, which is at the base of the analysis carried out in \cite{lamstr}.

 We note {\it en passant} that the standard eigenvalues of Maxwell's equations in a cavity are given by two families of positive numbers, one of which is the family of the squares of the zeros of the equation $j_l(k)=0$, see \cite{coda19} , or   \cite[Appendix]{lamzac},  for more details.
\end{remark}

\vspace{0.4cm}

\noindent
{\bf Acknowledgements}   The authors are thankful to the Departments of Mathematics of the University of Padova and of the National and Kapodistrian University of Athens for the kind hospitality.   In addition, they acknowledge financial support from the research project BIRD191739/19  ``Sensitivity analysis of partial differential equations in the mathematical theory of electromagnetism'' of the University of Padova.
The first named author (FF) is grateful for the received support to the UK Engineering and Physical Sciences Research Council through grant EP/T000902/1.
The second named author (PDL) is a  member of the Gruppo Nazionale per l'Analisi Matematica, la Probabilit\`a e le loro Applicazioni (GNAMPA) of the Istituto Nazionale di Alta Matematica (INdAM).

\end{document}